\documentclass{amsart}
\usepackage{graphicx,latexsym,amssymb,amsthm,amsmath}

 \newtheorem{thm}{Theorem}[section]
 \newtheorem{cor}[thm]{Corollary}
 \newtheorem{lem}[thm]{Lemma}
 \newtheorem{prop}[thm]{Proposition}
 \newtheorem{defn}[thm]{Definition}
 \newtheorem{rem}[thm]{Remark}
 \newtheorem{example}[thm]{Example}
 \newtheorem{question}[thm]{Question}

\newcommand{\ben}{\begin{enumerate}}
\newcommand{\een}{\end{enumerate}}
\newcommand{\ble}{\begin{lem}}
\newcommand{\ele}{\end{lem}}
\newcommand{\bth}{\begin{thm}}
\renewcommand{\eth}{\end{thm}}
\newcommand{\bpr}{\begin{prop}}
\newcommand{\epr}{\end{prop}}
\newcommand{\bco}{\begin{cor}}
\newcommand{\eco}{\end{cor}}
\newcommand{\bde}{\begin{defn}}
\newcommand{\ede}{\end{defn}}
\newcommand{\brem}{\begin{rem}}
\newcommand{\erem}{\end{rem}}
\newcommand{\bexm}{\begin{example}}
\newcommand{\eexm}{\end{example}}
\newcommand{\bqst}{\begin{question}}
\newcommand{\eqst}{\end{question}}

\newcommand{\C}{\mathbb{C}}
\newcommand{\D}{\mathbb{D}}
\newcommand{\PP}{\mathbb{P}}
\newcommand{\RR}{\mathbb{R}}

\begin{document}

\title{Several applications of Bezout matrices}

\author{Shmuel Kaplan, Alexander Shapiro and Mina Teicher}
\address{Shmuel Kaplan, Alexander Shapiro and Mina Teicher,
Department of Mathematics and Statistics, Bar-Ilan University,
Ramat-Gan 52900, Israel} \email{[kaplansh, sapial, teicher]@math.biu.ac.il}

\thanks{Partially supported by  EU-network HPRN-CT-2009-00099(EAGER), (The Emmy Noether Research Institute for Mathematics and the Minerva Foundation of Germany), the Israel Science Foundation grant \# 8008/02-3 (Excellency Center "Group  Theoretic Methods in the Study of Algebraic Varieties").}
\thanks{The authors wish to thank Prof. Victor Vinnikov for helpful advices.}

\maketitle

\begin{abstract}
The notion of Bezout matrix is an essential tool in studying broad variety of subjects:
zeroes of polynomials, stability of differential equations, rational transformations of
algebraic curves, systems of commuting nonselfadjoint operators, boundaries of quadrature
domains etc. We present a survey of several properties of Bezout matrices and their
applications in all mentioned topics. We use the framework of Vandermonde vectors because
such approach allows us to give new proofs of both classical and modern results and in
many cases to obtain new explicit formulas. These explicit formulas can significantly
simplify various computational problems and, in particular, make the research of
algebraic curves and their applications easier. In addition we wrote a \emph{Maple}
software package, which computes all the formulas. For instance, as Bezout matrices are
used in order to compute the image of a rational transformation of an algebraic curve, we
used these results to study some connections between small degree rational transformation
of an algebraic curve and the braid monodromy of its image.
\end{abstract}

\section*{Introduction}

Numerous works in operator theory shows that an algebraic curve given by a determinantal
representation can be associated to a system of commuting nonselfadjoint operators. In
the simplest case there are two commuting nonselfadjoint operators in a system and they
are rational images of the same operator. This particular case gives a motivation and a
framework for the studying the image of a complex line under a rational transformation.
Using the notion of a determinantal representation of an algebraic curve there was found
the explicit formula which describes the image of a complex line under a rational
transformation. The image is given by a determinantal representation which uses the
Bezout matrices of pairs of polynomials that define the rational transformation.

This property of Bezout matrices can be used for many purposes, i.e. to study the braid
monodromy of the image of two intersecting lines under rational transformations of small
degrees. Also using this property one can describe explicitly the boundary of a
quadrature domain. This property and many others can be proved using the Vandermonde
vectors, which are a natural framework to study Bezout matrices.

These properties were used by the authors of this paper to create a package of procedures
for \emph{Maple} software \cite{Ours2}. This package allows the computation of all
classical and modern formulas that are mentioned in the paper, and to perform many
different tasks.

\section{Bezout matrices and Vandermonde vectors}

First example of Bezout matrices for polynomials of small degree appeared in Euler's work in 1748,
\cite{Eul}. Using this example Bezout gave a general definition of Bezout matrices for polynomials of any
degree in 1764, \cite{Bez}. The notation of Bezoutian matrix was introduced by Sylvester in 1853,
\cite{Syl}. The most common definition was given by Cayley in 1857, \cite{Cay}.
\begin{defn}
For two polynomials in one variable $p(x)$ and $q(x)$ of degree $n$ there exists uniquely determined $n
\times n$ symmetric matrix $B(p,q) = \left( b_{ij} \right)_{i,j=1}^n$ such that
$$\frac{p(x)q(y) - q(x)p(y)}{x-y} = \sum_{i,j=1}^{n} b_{ij}x^{i-1}y^{j-1}$$
The matrix $B(p,q)$ is called the Bezout matrix of polynomials $p$ and $q$.
\end{defn}

We will use a framework of Vandermonde vectors. The matrices composed of Vandermonde vectors appeared
inexplicitly in Vandermonde's work in 1772, \cite{Van}. The notation was attributed to Vandermonde by Weill
in 1888, \cite{Wei}.
\begin{defn}
$$V_n(x) =
\left(
\begin{array}{c}
1 \\
x \\
\dots \\
x^{n-1}
\end{array}
\right) = \left( x^i \right)_{i=0}^{n-1}$$ $V_n(x)$ is called the Vandermonde vector of the length $n$.
\end{defn}

\begin{thm}
If $x_1, x_2, \dots , x_n$ are pairwise distinct then vectors \\
$V_n(x_1), V_n(x_2), \dots , V_n(x_n)$ are linearly independent.
\end{thm}
\begin{proof}
Let us consider $n \times n$ matrix
$$\left(V_n(x_1), V_n(x_2), \dots , V_n(x_n)\right) =
\left(
\begin{array}{cccc}
1 & 1 & \dots & 1 \\
x_1 & x_2 & \dots & x_n \\
\dots & \dots & \dots & \dots \\
x_1^{n-1} & x_2^{n-1} & \dots & x_n^{n-1}
\end{array}
\right)$$

Let us suppose that $p = (p_0, p_1, \dots , p_{n-1})$ is a row--vector from the left kernel of this matrix.
Then the polynomial $p(x) = p_0 + p_1 x + \dots +p_{n-1} x^{n-1}$ has $n$ zeroes. The degree of this
polynomial is less than $n$. Hence, the kernel of this matrix is trivial and vectors $V_n(x_1), V_n(x_2),
\dots , V_n(x_n)$ are linearly independent.
\end{proof}

We will call $V_n(x)$ the Vandermonde vector of order zero. It is natural to define the Vandermonde vector
of higher orders. For higher orders the definition is: $V^k_n(x) =
\left(\frac{d^k}{dx^k}x^i\right)_{i=0}^{n-1}$ We will call $V^k_n(x)$ the Vandermonde vector of order $k$.

\begin{thm}\label{thm:Linear independence of Vandermonde vectors}
If $i_1+i_2+\dots+i_m=n-1$ and $x_1, x_2, \dots , x_m$ are pairwise distinct then vectors $V_n(x_1),
V^1_n(x_1), \dots , V^{i_1}_n(x_1),
V_n(x_2), V^1_n(x_2), \dots , V^{i_2}_n(x_2),\\
\dots, V_n(x_m), V^1_n(x_m), \dots , V^{i_m}_n(x_m)$ are linearly independent.
\end{thm}
The prove is the same as of the above theorem.

\begin{lem}\label{lem:Basic properties}
For every two polynomials $p(x)$ and $q(x)$\\
1.$$V^T(x)B(p,q)V(y) = \frac{p(x)q(y) - q(x)p(y)}{x-y}$$ 2.
$$V^T(x)B(p,q)V(y)=V^T(x)(B(p,1)q(y)-B(q,1)p(y))V(y)$$ 3. $$V^T(x) B(p,q) V(x) = q(x)p'(x) - p(x)q'(x)$$
\end{lem}
\begin{proof} The first statement follows immediately from the definition of the Bezout
matrix. The second statement follows from the next decomposition:
$$\frac{p(x)q(y) - q(x)p(y)}{x-y} = \frac{p(x)q(y) - p(y)q(y) + p(y)q(y) - q(x)p(y)}{x-y} =$$
$$q(y)\frac{p(x) - p(y)}{x-y} - p(y)\frac{q(x) - q(y)}{x-y}=$$
$$q(y)V^T(x)B(p,1)V(y)-p(y)V^T(x)B(q,1)V(y)$$
To prove the third statement let us notice that
$$V^T(x) B(p,1) V(x) =\lim_{\epsilon \rightarrow 0} V^T(x) B(p,1) V(x+\epsilon) =$$
$$\lim_{\epsilon \rightarrow 0}\frac{p(x)-p(x+\epsilon)}{x-(x+\epsilon)}=
\lim_{\epsilon \rightarrow 0}\frac{p(x+\epsilon)-p(x)}{\epsilon}=p'(x)$$ This and the first statement of
the lemma imply the third statement.
\end{proof}
\begin{cor}
For every polynomials $p(x)$ and $q(x)$, vector $w$ and point $y$
$$w^T B(p,q) V(y) = w^T(B(p,1)q(y)-B(q,1)p(y))V(y)$$
\end{cor}
\begin{proof}
Vector $w$ can be represented as a linear combination of Vandermonde vectors. $w=\sum a_i V(z_i)$, where
all points $z_i$ are different from $y$. If follows from the second statement of the previous lemma that
$$w^T B(p,q) V(y) = \left(\sum a_i V^T(z_i)\right) B(p,q) V(y) = \sum a_i \left( V^T(z_i) B(p,q) V(y)
\right)=$$
$$\sum a_i \left( V^T(z_i)(B(p,1)q(y)-B(q,1)p(y))V(y)\right)=$$
$$\left( \sum a_i V^T(z_i)\right)(B(p,1)q(y)-B(q,1)p(y))V(y)=$$
$$w^T(B(p,1)q(y)-B(q,1)p(y))V(y)$$
\end{proof}

One of the very well known classical applications of Bezout matrices is Jacobi--Darboux theorem. This
theorem was proved independently by Jacobi in 1836, \cite{Jac} and by Darboux in 1876, \cite{Dar}.
\begin{thm}[Jacobi--Darboux theorem]
The number of common zeroes of two polynomials equals to the dimension of the kernel of the Bezout matrix
of these polynomials.
\end{thm}
There are different proofs of this theorem. For the proof in terms of Sylvester matrices and Vandermonde
vectors see \cite{ShV1}.

\section{Inverse of Bezout matrix and Hermite theorem}
In 1974 Lander proved that the inverse of Bezout matrix is a matrix of Hankel type, \cite{Lan}. This
property, its applications and similar properties of structured matrices were studied by Gohberg and
Olshevsky \cite{GoOl}, Heinig and Hellinger \cite{HeiHel}, Tyrtyshnikov \cite{Tyr}, Chen and Yang
\cite{ChYa}. The framework of Vandermonde vectors allows not only to reprove the result of Landau but also
to obtain explicit formulas for the coefficient of the inverse of Bezout matrix, \cite{ShV2}.

\begin{thm} \label{thm:Explicit formulas for Hankel matrices}
If two polynomials $p(x)$ and $q(x)$ of degree $n$ have no common zeroes and polynomial $p(x)$ has no
multiple zeroes then the inverse of Bezout matrix $B(p,q)$ is a Hankel type matrix $H=(h_{ij})$ and
$$h_{ij}=\sum_{k=1}^{n}\frac{x^{i+j-2}_{k}}{q(x_k)p'(x_k)},$$where $x_1, x_2 \dots x_n$ are zeroes of $p(x)$.
\end{thm}

\begin{proof}
Let us construct a matrix $V_{p}$ from Vandermonde vectors in zeroes of the polynomial $p(x)$:
$V_{p}=(V(x_1), V(x_2), \dots , V(x_n))$. We will denote by $D$ the matrix $V^{T}_{p}B(p,q)V_{p}$. It is
obvious that $B(p,q)=(V^{T}_{p})^{-1}D(V_{p})^{-1}$ and therefore $B^{-1}(p,q)=V_{p}D^{-1}V^{T}_{p}$.

Let us denote by $d_{ij}$ coefficients of the matrix $D$. If $i \neq j$ then
$$d_{ij}=V^{T}(x_i)B(p,q)V(x_j)=\frac{p(x_i)q(x_j)-q(x_i)p(x_j)}{x_i-x_j}=0,$$
since the numerator is zero and the denominator is not because all zeroes of $p(x)$ are different. If $i =
j$ then by the third statement of the Lemma \ref{lem:Basic properties}
$$d_{ii}=q(x_i)p'(x_i) - p(x_i)q'(x_i)=q(x_i)p'(x_i),$$
therefore $D$ is a diagonal matrix: $D=diag(q(x_1)p'(x_1), \dots , q(x_n)p'(x_n))$ which means that
$$D^{-1}=diag\left(\frac{1}{q(x_1)p'(x_1)}, \dots , \frac{1}{q(x_n)p'(x_n)}\right)$$ Therefore
$$V_{p}D^{-1}=\left(\frac{x^{i-1}_{k}}{q((x_k)p'(x_k)}\right)_{ik}$$
Clearly, $V^{T}_{p}=(x^{j-1}_{k})_{kj}$ Hence, the inverse of Bezout matrix is
$$B^{-1}(p,q)=V_{p}D^{-1}V^{T}_{p}=\left(\sum_{k=1}^{n}\frac{x^{i+j-2}_{k}}{q(x_k)p'(x_k)}\right)_{ij},$$
which proves the theorem.
\end{proof}

Hermite theorem was proved by Hermite in 1856, \cite{Her}. The theorem determines when all zeroes of a
polynomial belong to the upper--half plane. This theorem can be proved in a framework of Vandermonde
vectors, \cite{ShV2}.

\begin{thm}[Hermite theorem]
All zeroes of the polynomial $p(x)$ belongs to the upper--half plane if and only if the matrix
$\frac{1}{2i}B(p,\bar p)$  is positive definite.
\end{thm}

\section{Modern applications}

The simplest and very illustrative case of rational transformations of algebraic curves is a  of a rational
transformation of the complex line ${\C}$ into the complex plane ${\C}^2$. The image is a rational plane
algebraic curve and the explicit formulas for this image in terms of the polynomials that define the
rational transformation were obtained by Kravitsky in 1979, \cite{Kra}.
\begin{thm}
Three polynomials in two variables $p_0(x)$, $p_1(x)$ and $p_2(x)$ map complex line $\C$ into complex plane
$\C^2$:
$$x \to \left(\frac{p_1(x)}{p_0(x)},\frac{p_2(x)}{p_0(x)}\right)$$
The image is a rational curve defined by a polynomial
$${\Delta}(x_1,x_2)=\det(B(p_1,p_2)+x_1B(p_2,p_0)+x_2B(p_0,p_1))$$
\end{thm}
\begin{proof}
Let $x$ and $y$ be two points on the curve, $V(x)$ and $V(y)$ be two Vandermonde vectors.
$V^T(x)(p_0(x)B(p_1,p_2) + p_1(x)B(p_0,p_2) + p_2(x)B(p_0,p_1))V(y)=$\\
$p_0(x)V^T(x)B(p_1,p_2)V(y) + p_1(x)V^T(x)B(p_2,p_0)V(y) + p_2(x)V^T(x)B(p_0,p_1)V(y)=$
$(x-y)(p_0(x)(p_1(x)p_2(y)-p_2(x)p_1(y))+(p_2(x)p_0(y)-p_0(x)p_2(y))+(p_0(x)p_1(y)-p_1(x)p_0(y)))=0$ This
identity holds for for arbitrary $y$. Hence, $(p_0(x)B(p_1,p_2) + p_1(x)B(p_2,p_0) +
p_2(x)B(p_0,p_1))V(x)=0$ and therefore $$\det \left(B(p_1,p_2) + \frac{p_1(x)}{p_0(x)}B(p_2,p_0) +
\frac{p_2(x)}{p_0(x)}B(p_0,p_1)\right)=0,$$ which implies the theorem.
\end{proof}

\subsection{Nonselfadjoint operators}
The work of M. S. Liv\v{s}ic and his collaborators in operator theory associates to a system of commuting
nonselfadjoint operators an algebraic curve (called the discriminant curve) given by a determinantal
representation, see \cite{LMKV}. This discovery leads to a very fruitful interplay between operator theory
and algebraic geometry: problems of operator theory lead to problems of algebraic geometry and vice versa.

A natural problem in operator theory is to define properly the notion of a rational transformation of a
system of commuting nonselfadjoint operators. This arises whenever one wants to study the algebra generated
by a given system of commuting nonselfadjoint operators. It may also allow representing the given system of
commuting nonselfadjoint operators in terms of another system which is simpler in some sense (e.g., it
contains fewer operators, or the operators have a smaller nonhermitian rank). A related problem in
algebraic geometry is to find an image of an algebraic curve given by a determinantal representation under
a rational transformation.

To formulate the main problem more precisely, we have to introduce some notation. We
shall use the framework of commutative vessels \cite{LMKV} which turns out to be very
convenient in the study of commuting nonselfadjoint operators; it generalizes to the
multi--operator case the framework of colligations (nodes) which has been extensively
used in the study of a single nonselfadjoint (or non-unitary) operator, see, e.g.,
\cite{Bro}.

Let $H$ be a Hilbert space (finite-- or infinite--dimensional) and let $E$ be a finite--dimensional Hilbert
space.

\begin{defn}
An operator node is a collection
\begin{eqnarray}
C=(A, H, {\Phi}, E, {\sigma}) \nonumber
\end{eqnarray}
where $A:H \rightarrow H$ is bounded linear operator, ${\Phi}$ is a bounded linear mapping from $H$ to $E$
with the adjoint mapping ${\Phi}^{*}:E \rightarrow H$, ${\sigma}$ is bounded selfadjoint operators in $E$,
such that ${\Phi}^*{\sigma}{\Phi} = \frac{ 1 }{ i }(A-A^*)$.
\end{defn}

\begin{defn} A commutative vessel is a collection
\begin{eqnarray}
V=(A_1, A_2, H, {\Phi}, E, {\sigma}_1, {\sigma}_2, {\gamma}^{in}, {\gamma}^{out}) \nonumber
\end{eqnarray}
where $A_1, A_2: H \rightarrow H$ are bounded linear commuting operators, ${\Phi}$ is a bounded linear
mapping from $H$ to $E$ with the adjoint mapping ${\Phi}^{*}:E \rightarrow H$, ${\sigma}_1$, ${\sigma}_2$,
${\gamma}^{in}$, ${\gamma}^{out}$ are bounded selfadjoint operators in $E$, such that
${\gamma}^{in}=-{\gamma}^{in}, {\gamma}^{out}=-{\gamma}^{out}$, and
\begin{eqnarray}
{\Phi}^*{\sigma}_{i}{\Phi} & = & \frac{ 1 }{ i }(A_{i}-A^*_{i}) \nonumber \\
{\gamma}^{in}{\Phi} & = & {\sigma}_1{\Phi}A^*_{2}-{\sigma}_{2}{\Phi}A^*_{1} \nonumber \\
{\gamma}^{out}{\Phi} & = & {\sigma}_1{\Phi}A_{2}-{\sigma}_{2}{\Phi}A_{1}\nonumber \\
{\gamma}^{out} & = & {\gamma}^{in}+i({\sigma}_{1}{\Phi}{\Phi}^*{\sigma}_{2}-
{\sigma}_{2}{\Phi}{\Phi}^*{\sigma}_{1}) \nonumber
\end{eqnarray}
\end{defn}

For the simplest case, when a system of operators consists of a single operator and the discriminant curve
is a line, these problems were solved by Kravitsky in 1979, \cite{Kra}. Let us consider an operator node
$C=(A, H, {\Phi}, E, {\sigma})$ and three polynomials $p_0(x)$, $p_1(x)$ and $p_2(x)$ of degree $n$, such
that $p_0(A)$ is an invertible operator. We define $A_1=p_1(A)p_0^{-1}(A)$ and $A_2=p_2(A)p_0^{-1}(A)$
\begin{thm}
A collection $V=(A_1, A_2, H, {\Phi}', E', {\sigma}_1, {\sigma}_2, {\gamma}^{in}, {\gamma}^{out})$, where
\begin{eqnarray}
E' & = & E^{\otimes n} \nonumber \\
{\Phi}' & = & P_{E}({\Phi}A^{*i_1}_1A^{*i_2}_2p^{-1}_0(A^*_1A^*_2))_{i_1i_2} \nonumber \\
{\sigma}_{1} & = & B(p_0,p_1)\otimes{\sigma} \nonumber \\
{\sigma}_{2} & = & B(p_0,p_2)\otimes{\sigma} \nonumber \\
{\gamma}^{in} & = &  B(p_1,p_2)\otimes{\sigma} \nonumber \\
{\gamma}^{out} & = &
{\gamma}^{in}+i({\sigma}_1{\Phi}'{\Phi}'^*{\sigma}_2-{\sigma}_2{\Phi}'{\Phi}'^*{\sigma}_1) \nonumber
\end{eqnarray}
is a vessel and the the discriminant curve of this vessel is defined by the equation
$\det^{m}(B(p_1,p_2)+x_1B(p_2,p_0)+x_2B(p_0,p_1))=0$
\end{thm}

Explicit formulas for an image of a plane algebraic curve given by a determinantal representation under a
rational transformation were obtained using generalization of Bezout matrices in \cite{ShV1}. These results
allowed to define properly the rational image of a system of two operators where the discriminant curve is
a plane algebraic curve, \cite{ShV3}.

\subsection{Quadrature domains}
The last two decades have witnessed a renewed interest and constant progress in the theory of quadrature
domains. Questions such as the constructions and parametrization of quadrature domains with prescribed
distribution u, the algebraic structure of the boundary or various functional analytic characterizations of
quadrature domains have been successfully investigated, see \cite{Sha}, \cite{Put}. One of possible ways to
find a polynomial that defines the boundary of a quadrature domain is to consider such rational
transformations of a complex plane that the image is a Riemann surface equipped with an involution.

The domain ${\Omega} \subset {\C}$ is called a quadrature domain if there exists a distribution $u$ with
finite support in ${\Omega}$ such that
$$\int_{\Omega}fdA = u(f),$$
for every integrable analytic function $f$ in ${\Omega}$. To be more specific, there are points
${\lambda}_j\in{\Omega}$ and constants ${\gamma}_{jk}$, $0 \leq k \leq m(j)-1$, $1 \leq j \leq m$ such that
$$\int_{\Omega}fdA = \sum_{j=1}^m \sum_{k=0}^{m(j)-1}{\gamma}_{jk}f^{(k)}({\lambda}_j),$$
where $m(j)$ is the multiplicity of the point ${\lambda}_j$. The quadrature domain is an image of the unit
disk ${\D}$ under polynomial transformation defined by a polynomial $q(z)$ and the boundary of the
quadrature domain is given by the equation ${\Delta}(z, \overline{z})=0$. To find the polynomial ${\Delta}$
we consider a rational transformation of the complex plane: $z \rightarrow \left(q(z),
\overline{q(\frac{1}{\overline{z}})}\right)$.
\begin{thm}
If a quadrature domain ${\Omega}$ is an image of the unit disk ${\D}$ under polynomial transformation
defined by a polynomial $q(z)$ then there exist three polynomials $p_0(z)$, $p_1(z)$ and $p_2(z)$ such that
$q(z)=\frac{p_1(z)}{p_0(z)}$ and $\overline{q(\frac{1}{\overline{z}})} = \frac{p_2(z)}{p_0(z)}$. The
boundary of ${\Omega}$ is defined by the equation $\det(B(p_1,p_2)+zB(p_2,p_0)+\overline{z}B(p_0,p_1))=0$.
\end{thm}
For the proof see \cite{PShV}.

\subsection{Braid monodromy}

In 1937 Zariski laid down the foundations for the braid monodromy of curves in $\C\PP^2$. The braid
monodromy is a homomorphism between the fundamental group of a punctured disk and the braid group. We give
here a short description of the braid monodromy.

Let us recall that the braid group $B_n$ is the $\mathcal{MCG}$ (mapping class group) of the $n$-puncture
disk. We distinguish some important elements in the braid group $B_n$ which are called \emph{half-twists}.
Let $D$ be a closed disk and let $K=\{k_1,\cdots,k_n\} \subset D \setminus \partial D$. Choose $u\in
\partial D$. Let $a,b$ be two points of $K$. We denote $K_{a,b}=K \setminus \{a,b\}$. Let $\sigma$ be a
simple path in $D \setminus (\partial D \cup K_{a,b})$ connecting $a$ with $b$. Choose a small regular
neighborhood $U$ of $\sigma$ and an orientation
preserving diffeomorphism $f:\RR ^2 \to \C$ such that $f(\sigma )=[-1,1]$, $f(U)=\{z \in \C \ | \ 
|z| <2\}$.

Let $\alpha (x)$, $0 \leq x$ be a real smooth monotone function such that:

$\alpha (x)=\left \{
\begin{array}{ll}
1, & 0 \leq x \leq \frac{3}{2} \\
0, & 2 \leq x
\end{array} \right.$

Define a diffeomorphism $h:\C \to \C$ as follows: for $z=re^{i\varphi} \in \C$ let
$h(z)=re^{i(\varphi +\alpha (r)\pi )}$

For the set $\{z \in \C \ | \ 2 \leq |z|\}$,  $h(z)={\rm Id}$,
and for the set $\{z \in \C \ | \ |z|\leq \frac{3}{2}\}$, $h(z)$ is a rotation by
$180 ^{\circ}$ in the positive direction.

Considering $(f \circ h \circ f^{-1})|_D$ (we will compose from left to
right) we get a diffeomorphism of $D$ which switches $a$ and $b$ and is the identity on $D
\setminus U$. Thus it defines an element of $B_n[D,K]$.

The diffeomorphism $(f \circ h \circ f^{-1})|_D$ defined above induces an automorphism on $\pi _1(D
\setminus K,u)$, that switches the position of two generators of $\pi _1(D \setminus K,u)$.

\bde
Let $H(\sigma)$ be the braid defined by $(f \circ h \circ f^{-1})|_D$. We call
$H(\sigma )$ the \emph{positive half-twist defined by $\sigma$}.
\ede

Let $C$ be a real curve in $\C^2$ of degree $n$. Denote by $pr_1:C \to \C$ and by $pr_2:C \to \C$ the
projections to the first and second coordinate, defined in the obvious way. For $x \in \C$ we denote $K(x)$
the projection of the points in $C$ which lie with $x$ as their first coordinate to the second coordinate
(i.e., $K(x)=pr _2(pr _1^{-1}(x))$).

Let $N \subset \C$ be the set $N=N(C)=\{x \in\C \ | \ |K(x)|<n\}=\{x_1, \cdots ,x_p\}$. We restrict
ourselves only to the cases where $N$ is finite. Take $E$ to be a closed disc in $\C$ for which $N \subset
E \setminus \partial E$. In addition take $D$ to be a closed disc in $\C$ for which $D$ contains all the
points $\{K(x) \ | \ x\in E\}$. That means that when restricted to $E$, we have $C \subset E \times D$.

With these definitions in hand we may define the \emph{braid monodromy of a projective curve}:

\begin{defn} \label{def:braid monodromy}
Let $C$ be a projective curve of degree $n$ in $\C\PP^2$, $L$ be a generic line at infinity such that $|L
\cap C|=n$, and $(x,y)$ is an affine coordinate system for $\C^2=\C\PP^2 \setminus L$ such that the
projection of $C$ to the first coordinate is generic. For $E,D,N$ defined as above, let $M \in \partial E
\cap \RR$ be the base point of $\pi _1(E \setminus N)$, and let $\sigma$ be an element of $\pi _1(E
\setminus N)$. To $\sigma$ there are $n$ lifts in $C$, each one of them begins and ends in the points of $M
\times K(M)$. Projecting these lifts using $pr _2:C \to \C$ we get $n$ paths in $D$ which begin and end in
the points of $K(M)$. These induce a diffeomorphism of $\pi _1(D \setminus K(M))$ which is the braid group
$B_n$ as defined earlier. We call the homomorphism $\varphi :\pi _1(E \setminus N) \to B_n$ the \emph{braid
monodromy} of $C$ with respect to $L,E \times D,pr _1,$ and $M$.
\end{defn}

It is natural to ask questions about the connection between the rational transformation and the braid
monodromy induced by its image. For example one can formulate the following questions:

\bqst Study the singular points of the image of a rational transformation, and define conditions on
$p_0(x,y),p_1(x,y),p_2(x,y)$ which induce specific braid monodromy results. \eqst

\bqst
Let $C$ be a curve. let $r$ be a rational transformation. Classify all braid monodromy results which may result.
\eqst

\bqst
Formulate necessary and sufficient conditions for a braid monodromy to be of rational curve.
\eqst

\bqst Given a braid monodromy, which satisfies the sufficient condition above. Formulate a family of
rational curves which will induce such braid monodromy. \eqst

\bqst
Given two isomorphic rational curves. What can be said on their braid monodromies.
\eqst

Of course, these questions are more than wide, and at this point may not be completely answered. In
\cite{Ours} we established two results concerning degree 2 rational transformations, as follows:

Let us consider the rational transformation
$$(x,y) \mapsto (p_0(x,y),p_1(x,y),p_2(x,y)).$$
In order to compute the local braid monodromy at the point $(p_0(x_0,y_0)$, $p_1(x_0,y_0)$, $p_2(x_0,y_0))$,
we assume that $p_0(x_0,y_0) \neq 0$.
We define:
$r_1(x,y)=\frac{p_1(x,y)}{p_0(x,y)}$, $r_2(x,y)=\frac{p_2(x,y)}{p_0(x,y)}$, and recursively \\
$D_1(x)=r_2'(x,0) \cdot \frac{1}{r_1'(x,0)}$\\
$D_n(x)=D_{n-1}'(x,0) \cdot \frac{1}{r_1'(x,0)}$\\
$E_1(y)=r_2'(0,y) \cdot \frac{1}{r_1'(0,y)}$\\
$E_n(y)=D_{n-1}'(0,y) \cdot \frac{1}{r_1'(0,y)}$\\

\begin{cor} \cite{Ours} \label{cor:BM_Criteria}
Let $(p_0(x_0,y_0),p_1(x_0,y_0),p_2(x_0,y_0))$ be one of the intersection points of the two conics at the
image $r(C)$. Let $i$ be the minimal index for which $D_i(x_0) \neq E_i(y_0)$. Then, the multiplicity of
the intersection point is $i+1$, and thus the local braid monodromy at this intersection point is $(i+1)$
full twists of two strings.
\end{cor}

\begin{thm} \cite{Ours} \label{thm:global braid moodromy}
Let $C$ be a curve which consists of two intersecting lines, and let $r$ be a real rational transformation
of degree $2$. Then, the braid monodromy of $r(C)$ is completely defined by the number and multiplicity of
it's real self intersection points.
\end{thm}

Theorem \ref{thm:global braid moodromy} gives a full classification of the braid monodromy of the image of
two intersecting lines under degree 2 rational transformations.

\end{document}